\documentclass[12pt]{amsart}
\usepackage{times,fullpage}%,backref}
\usepackage{latexsym,amscd,amssymb}
\newtheorem{thm}{Theorem}

\theoremstyle{remark}

\theoremstyle{definition}

\newcommand{\C}{\mathbb{ C}}

\newcommand{\Q}{\mathbb{ Q}}

\newcommand{\CC}{\mathcal{ C}}

\newcommand{\PP}{\mathcal{EP}}
\newcommand{\TT}{\mathcal{HT}}

\begin{document}

\title{Characteristic numbers of algebraic varieties}
\author{D.~Kotschick}
\address{Mathematisches Institut, {\smaller LMU} M\"unchen,
Theresienstr.~39, 80333~M\"unchen, Germany}
\email{dieter@member.ams.org}
\date{June 23, 2009; published in Proc.~Natl.~Acad.~Sci.~USA {\bf 106} (2009), no. 25, 10114--10115. }

\maketitle

\begin{abstract} 
A rational linear combination of Chern numbers is an oriented diffeomorphism invariant of smooth complex projective 
varieties if and only if it is a linear combination of the Euler and Pontryagin numbers. In dimension at least three
only multiples of the top Chern number, which is the Euler characteristic, are invariant under diffeomorphisms 
that are not necessarily orientation-preserving. 

In the space of Chern numbers there are two distinguished subspaces, one spanned by the Euler and Pontryagin numbers, 
the other spanned by the Hirzebruch--Todd numbers. Their intersection is the span of the Euler number and the signature. 
 \end{abstract}

\bigskip

\section{Introduction}

In 1954, Hirzebruch~\cite[Problem~31]{Hir1} asked which linear combinations of Chern numbers of smooth
complex projective varieties are topologically invariant. 
The purpose of this article is to announce a complete solution to this problem, and to discuss some additional
results about the characteristic numbers of algebraic varieties that arise from this solution. 
Some of these further results characterize the characteristic numbers of smooth complex projective varieties that
can be bounded in terms of Betti numbers. Another one of these results, Theorem~\ref{t:main} below, 
determines the intersection of Pontryagin and Hirzebruch--Todd numbers in even complex dimensions.

\section{Solution of Hirzebruch's problem}

Because the manifold underlying a complex-algebraic variety has a preferred orientation, it is most natural to 
interpret the term ``topological invariance'' in Hirzebruch's problem to mean invariance under orientation-preserving 
homeo- or diffeomorphisms. 
With this interpretation the solution to the problem is given by:
\begin{thm}\label{t:DO}
A rational linear combination of Chern numbers is an oriented diffeomorphism invariant of smooth complex projective 
varieties if and only if it is a linear combination of the Euler and Pontryagin numbers. 
\end{thm}
In one direction, the Euler number, which is the top Chern number, is of course a homotopy invariant.
Further, the Pontryagin numbers, which are special linear combinations of Chern numbers, are oriented
diffeomorphism invariants. In fact, Novikov proved that the Pontryagin numbers are also invariant under 
orientation-preserving homeomorphisms, and so Theorem~\ref{t:DO} is unchanged if we replace 
oriented diffeomorphism-invariance by oriented homeomorphism-invariance.
The other direction, proving that there are no other linear combinations that are oriented diffeo\-mor\-phism-invariants,
has proved to be quite difficult because of the scarcity of examples of diffeomorphic algebraic
varieties with distinct Chern numbers.

Given Theorem~\ref{t:DO}, and the fact that Pontryagin numbers depend on the orientation, one might expect that
only the Euler number is invariant under homeo- or diffeomorphisms that do not necessarily preserve the 
orientation. For diffeomorphisms this is almost but not quite true:
\begin{thm}\label{t:D}
In complex dimension $n\geq 3$ a rational linear combination of Chern numbers is a diffeomorphism invariant 
of smooth complex projective varieties if and only if it is a multiple of the Euler number $c_n$. In complex dimension
$2$ both Chern numbers $c_2$ and $c_1^2$ are diffeomorphism invariants of complex surfaces.
\end{thm}
The statement about complex dimension $2$ is a consequence of Seiberg--Witten theory.
It is an exception due to the special nature of 
four-dimensional differential topology. The exception disappears if we consider homeomorphisms instead of diffeomorphisms:
\begin{thm}\label{t:H}
A rational linear combination of Chern numbers is a homeomorphism invariant 
of smooth complex projective varieties if and only if it is a multiple of the Euler number. 
\end{thm}

These theorems show that linear combinations of Chern numbers of complex projective varieties are not usually determined
by the underlying manifold. This motivates the investigation of a modification of Hirzebruch's original problem, asking how far 
this indeterminacy goes. More precisely, one would like to know which linear combinations of Chern numbers are determined 
up to finite ambiguity by the topology. The obvious examples for which this is true are the numbers
\begin{equation}\label{eq:HRR}
\chi_{p} = \chi (\Omega^p)= \sum_{q=0}^{n} (-1)^{q}h^{p,q} \ .
\end{equation}
By the Hirzebruch--Riemann--Roch theorem~\cite{TMAG} these are indeed linear combinations of Chern numbers.
By their very definition, together with the Hodge decomposition of the cohomology, the numbers $\chi_p$ are bounded
above and below by linear combinations of Betti numbers. It turns out that this property characterizes the linear 
combinations of the $\chi_p$, as shown by the following theorem:
\begin{thm}\label{t:oBetti}
A rational linear combination of Chern numbers of complex projective varieties can be bounded in terms of Betti 
numbers if and only if it is a linear combination of the $\chi_p$. 
\end{thm}

%As a consequence of this result, most linear combinations of Chern numbers are independent of the Hodge structure:
%\begin{cor}\label{c:Hodge}
%A rational linear combination of Chern numbers of complex-algebraic varieties is determined by the Hodge  
%numbers if and only if it is a linear combination of the $\chi_p \ $. 
%\end{cor}
The span of the $\chi_p$ includes the Euler number $c_n = \sum_p (-1)^p\chi_p$
and the signature, which, according to the Hodge index theorem, equals $\sum_p \chi_p$.
It also includes the Chern number $c_1c_{n-1}$, by a result of Libgober and Wood~\cite[Theorem~3]{LW}.
Nevertheless, the span of the $\chi_p$ is a very small subspace of the space of linear combinations of Chern
numbers. The latter has dimension equal to $\pi(n)$, the number of partitions of $n$, which grows exponentially 
with $n$. The former has dimension $[(n+2)/2]$, the integral part of $(n+2)/2$. This follows from the symmetries 
of the Hodge decomposition, which imply $\chi_p = (-1)^n\chi_{n-p}$, together with the linear independence of
$\chi_0,\ldots,\chi_{[n/2]}$, which can be checked, for example, by evaluating on all $n$-dimensional products 
of $\C P^1$ and $\C P^2$.

The proofs of Theorems~\ref{t:DO}, \ref{t:D}, \ref{t:H} and~\ref{t:oBetti} are carried out in~\cite{CC}.

\section{Euler--Pontryagin versus Hirzebruch--Todd}

We now restrict ourselves to even complex dimensions $n=2m$.
Consider the $\Q$-vector space $\CC_{2m}$ of linear combinations of Chern numbers of almost complex manifolds of complex dimension $2m$, 
equivalently of real dimension $4m$. This is a vector space of dimension $\pi (2m)$, the number of partitions of $2m$, with a basis given by the monomials
$c_{i_1}\cdot\ldots\cdot c_{i_k}$ with $\sum_j i_j=2m$.

Define $\PP_{2m}\subset\CC_{2m}$ to be the subspace spanned by the Euler number $e=c_{2m}$ and by the Pontryagin numbers 
$p_{i_1}\cdot\ldots\cdot p_{i_k}$ with $\sum_j i_j=m$, where
$$
p_i = c_i^2-2c_{i-1}c_{i+1}+2c_{i-2}c_{i+2}+\cdots+(-1)^i 2c_{2i}
$$
are the Pontryagin classes. This subspace has dimension $1+\pi (m)$.

In the same spirit as Theorem~\ref{t:oBetti} we have the following:
\begin{thm}\label{t:Betti}
A rational linear combination of Pontryagin numbers of complex projective varieties can be bounded in terms of Betti 
numbers if and only if it is a multiple of the signature. 
\end{thm}
In other words, an element of $\PP_{2m}$ evaluated on projective varieties can be bounded in terms of Betti numbers if and only if it is a linear combination of the Euler 
number and the signature.

Theorem~\ref{t:Betti} is implicit in the papers~\cite{CC,NNP}, but is not explicitly stated there. Corollary~3 of~\cite{NNP} says that a
rational linear combination of Pontryagin numbers of smooth manifolds can be bounded in terms of Betti numbers if and only if it is a 
multiple of the signature. The proof given in~\cite{NNP} uses sequences of ring generators for the oriented bordism ring
 $\Omega_{\star}\otimes\Q$ which, in dimensions $\geq 8$, are $\C P^{2k}$-bundles over $S^4$, and are not 
complex-algebraic varieties. However, in that proof 
one can replace $S^4$ by $T^4$ endowed with the algebraic structure of an Abelian surface, and one can assume that the 
$\C P^{2k}$-bundles over it are algebraic, so as to obtain a proof of Theorem~\ref{t:Betti} above. This amounts to mapping 
the proof of Theorem~\ref{t:oBetti} given in~\cite{CC} from $\Omega_{\star}^U\otimes\Q$ to $\Omega_{\star}\otimes\Q$ under the forgetful natural
transformation.

Next define $\TT_{2m}\subset\CC_{2m}$ to be the subspace spanned by the Hirzebruch--Todd numbers $T^p_{2m}$ defined in~\cite[1.8]{TMAG}.
These are certain universal linear combinations of Chern numbers, which occur in the expansion
\begin{equation}\label{eq:T}
T_{2m}(y;c_1,\ldots,c_{2m}) = \sum_{p=0}^{2m}T^p_{2m}(c_1,\ldots,c_{2m})y^p \ ,
\end{equation}
where $T_{2m}(y;c_1,\ldots,c_{2m})$ is the multiplicative sequence corresponding to the power series
$$
Q(y;x) = x +\frac{x(y+1)}{e^{x(y+1)}-1} \ .
$$
It follows from the definition that $T^p_{2m}=T^{2m-p}_{2m}$, so that there are at most $m+1$ distinct $T^p_{2m}$, namely $T^0_{2m},\ldots,T^m_{2m}$.
By evaluating on products of projective spaces one sees that these $m+1$ linear combinations of Chern numbers are indeed linearly independent over $\Q$, 
so that $\TT_{2m}$ is of dimension $m+1$.

We have the following result about the relation between Pontryagin and Hirzebruch--Todd numbers:
\begin{thm}\label{t:main}
The intersection $\PP_{2m}\cap\TT_{2m}$ is two-dimensional, spanned by the Euler number and the signature.
\end{thm}

The statement is vacuous for $m=1$, as $\CC_2=\PP_2=\TT_2$.
In all dimensions, it is clear that the Euler number $c_{2m}=\sum_p(-1)^p T^p_{2m}$ is in the intersection. 
The signature is a linear combination of Pontryagin numbers by a result of Thom, and so is 
in $\PP_{2m}$. The space $\TT_{2m}$ of Hirzebruch--Todd numbers contains the $L$-genus $L_m$; cf.~\cite[1.8]{TMAG}. By the Hirzebruch
signature theorem $L_m$ equals the signature. So to prove Theorem~\ref{t:main} we have to prove that the only Pontryagin numbers in 
$\PP_{2m}\cap\TT_{2m}$ are the multiples of the signature.

\begin{proof}[The first proof of Theorem~\ref{t:main}.]
The Hirzebruch--Riemann--Roch theorem~\cite{TMAG} implies that in complex dimension $2m$ the $\chi_p$-numbers in~\eqref{eq:HRR}
are equal to the Hirzebruch--Todd numbers $T^p_{2m}$. In the special case when we evaluate on a smooth complex projective variety, 
or, more generally, on a compact K\"ahler manifold, the Hodge numbers $h^{p,q}$ are bounded by the Betti numbers
because of the Hodge decomposition of cohomology. Therefore, in this case, the Hirzebruch--Todd numbers $T^p_{2m}$ are 
bounded by linear combinations of Betti numbers, and Theorem~\ref{t:main} follows from Theorem~\ref{t:Betti}.
This completes our first proof of Theorem~\ref{t:main}.
\end{proof}

One can also give a more direct proof of Theorem~\ref{t:main}, without appealing to Theorem~\ref{t:Betti}. 

\begin{proof}[The second proof of Theorem~\ref{t:main}.]
The rational complex cobordism ring $\Omega_{\star}^U\otimes\Q$ is a polynomial ring on generators $\gamma_1,\gamma_2,\cdots$. The 
proof of Theorem~\ref{t:oBetti} in~\cite{CC} provides a specific choice for the $\gamma_i$, which, for $i\geq 3$, are in the kernel of the 
Hirzebruch $\chi_y$-genus defined by $\chi_y = \sum_{p} \chi_p y^p$.
In fact, since products of $\C P^1$ and $\C P^2$ detect all the components of the $\chi_y$-genus, any basis sequence can be modified by 
the addition of such products to achieve $\gamma_i\in\ker\chi_y$ for $i\geq 3$.

The Hirzebruch--Riemann--Roch theorem~\cite{TMAG} implies that the $\chi_p$-numbers agree with the Hir\-ze\-bruch--Todd numbers $T^p_{2m}$,
so that the subspace $\ker\chi_y\subset \Omega_{2m}^U\otimes\Q$, which is spanned by monomials containing a $\gamma_i$ with $i\geq 3$, 
coincides with the annihilator of $\TT_{2m}$.

The annihilator of the Pontryagin numbers is the kernel of the forgetful map
$\Omega_{\star}^U\otimes\Q\longrightarrow\Omega_{\star}\otimes\Q$,
which is the ideal generated by the $\gamma_i$ with odd $i$.

Now, if $f\in \PP_{2m}\cap\TT_{2m}$, then its kernel in the bordism group $\Omega_{2m}^U\otimes\Q$ contains all monomials containing a 
$\gamma_i$ with $i\geq 3$, and it also contains all monomials containing a $\gamma_i$ with $i$ odd. Therefore, the only monomial in the 
generators not contained in $\ker f$ is $\gamma_2^m$. If we normalize $f$, then it follows that $f=L_m$.
This completes our second proof of Theorem~\ref{t:main}.
\end{proof}

It would be interesting to know whether there is an elementary algebraic or combinatorial proof of Theorem~\ref{t:main}, showing directly from the definition~\eqref{eq:T} of 
the Hirzebruch--Todd numbers that their span contains no Pontryagin numbers other than the multiples of $L_m$, without using bordism theory or the 
Hirzebruch--Riemann--Roch theorem.

\subsection*{Acknowledgments}
I gratefully acknowledge several helpful comments and suggestions from P.~Landweber and the support of The Bell Companies Fellowship at the Institute for Advanced Study in Princeton.

\end{document}